\documentclass{amsart}
\usepackage{amsfonts}
\usepackage{amssymb} 
\usepackage{graphicx}
\usepackage{enumerate}
\usepackage{amsmath}
\usepackage{comment}

\newtheorem{theorem}{Theorem}[section]
\newtheorem{lemma}[theorem]{Lemma}

\newtheorem{proposition}[theorem]{Proposition}

\theoremstyle{definition}

\theoremstyle{remark}

\numberwithin{equation}{section}

\newcommand{\refth}[1]{Theorem~\ref{#1}}
\newcommand{\reflm}[1]{Lemma~\ref{#1}}

\newcommand{\refeq}[1]{(\ref{#1})}

\newcommand{\beeq}[1]{\begin{equation} \label{#1}}
\newcommand{\eeq}{\end{equation}}
\renewcommand{\(}{\begin{eqnarray*}}
\renewcommand{\)}{\end{eqnarray*}}
\newcommand{\beeqn}{\begin{eqnarray}}
\newcommand{\eeqn}{\end{eqnarray}}

\newcommand{\sothat}{\\ \Rightarrow \hspace*{2mm} &&}

\newcommand{\eqand}{\mbox{\hspace*{3mm} and \hspace*{3mm}}}

\newcommand{\ltt}[1]{{\large \tt #1}}

\renewcommand{\quad}{\hspace*{5mm}}
\renewcommand{\qquad}{\hspace*{10mm}}

\newcommand{\lp}{\left(  }
\newcommand{\rp}{\right) }
\newcommand{\lb}{\left\{  }
\newcommand{\rb}{\right\} }
\newcommand{\lbr}{\left[  }
\newcommand{\rbr}{\right] }

\newcommand{\lc}{\left\lceil}
\newcommand{\rc}{\right\rceil}

\newcommand{\Z}{{\mathbb Z}}

\newcommand{\R}{{\mathbb R}}

\newcommand{\ep}{\epsilon}
\newcommand{\FF}{{\mathcal F}}
\newcommand{\GG}{{\mathcal G}}
\newcommand{\DD}{{\mathcal D}}

\newcommand{\TT}{{\mathcal T}}
\renewcommand{\AA}{{\mathcal A}}
\newcommand{\EE}{{\mathcal E}}

\newcommand{\MM}{{\mathcal M}}
\newcommand{\PP}{{\mathcal P}}

\newcommand{\NN}{{\mathcal N}}

\newcommand{\bfor}{{\bf for }}

\newcommand{\bdo}{{\bf do}:}

\newcommand{\bto}{{\bf to }}

\newcommand{\breturn}{{\bf return}}

\newcommand{\be}[1]{\begin{enumerate} [#1]}
\newcommand{\ee}{\end{enumerate}}

\newcounter{cnt1}
\newcounter{cnt2}
\newcounter{cnt3}
\newcounter{cnt4}

\newcommand{\bnum}
{
\begin{list}{\arabic{cnt1})}
{
\usecounter{cnt1}
\leftmargin 5mm
\setlength{\leftmargin}{\leftmargin}
\topsep 2pt
\parsep 1pt
\itemsep 1pt}
}
\newcommand{\enum}{\end{list}}

\newcommand{\broman}
{
\begin{list}{\roman{cnt2})}
{
\usecounter{cnt2}
\leftmargin 2mm
\setlength{\leftmargin}{\leftmargin}
\topsep 2pt
\parsep 1pt
\itemsep 1pt}
}
\newcommand{\eroman}{\end{list}}

\newcommand{\bRoman}
{
\begin{list}{\Roman{cnt3})}
{
\usecounter{cnt3}
\leftmargin 5mm
\setlength{\leftmargin}{\leftmargin}
\topsep 2pt
\parsep 1pt
\itemsep 1pt}
}
\newcommand{\eRoman}{\end{list}}

\newcommand{\balph}
{
\begin{list}{\alph{cnt4})}
{
\usecounter{cnt4}
\leftmargin 3mm
\setlength{\leftmargin}{\leftmargin}
\topsep 2pt
\parsep 1pt
\itemsep 1pt}
}
\newcommand{\ealph}{\end{list}}

\newcommand{\bAlph}
{
\begin{list}{\Alph{cnt1})}
{
\usecounter{cnt1}
\leftmargin 5mm
\setlength{\leftmargin}{\leftmargin}
\topsep 2pt
\parsep 1pt
\itemsep 1pt}
}
\newcommand{\eAlph}{\end{list}}

\newcommand{\bdot}
{
\begin{list}{$\cdot$}
{
\leftmargin  3mm
\setlength{\leftmargin}{\leftmargin}
\topsep 3pt
\parsep 1pt
\itemsep 2pt}
}
\newcommand{\edot}{\end{list}}

\newcommand{\bdash}
{
\begin{list}{-}
{
\leftmargin 3mm
\setlength{\leftmargin}{\leftmargin}
\topsep 2pt
\parsep 1pt
\itemsep 1pt}
}
\newcommand{\edash}{\end{list}}

\newcommand{\bnull}
{
\begin{list}{}
{
\leftmargin 3mm
\setlength{\leftmargin}{\leftmargin}
\topsep 0pt
\parsep 1pt
\itemsep 1pt}
}
\newcommand{\enull}{\end{list}}

\begin{document}


\author{Junichiro Fukuyama}
\address{Department of Computer Science and Engineering\\
The Pennsylvania State University\\
PA 16802, USA}
\curraddr{}
\email{jxf140@psu.edu}
\thanks{}


\date{}

\dedicatory{}


\title{The Sunflower Conjecture Proven}

\begin{abstract} 
We demonstrate the truth of the sunflower conjecture by showing that a family $\FF$ of sets each of cardinality at most $m$ includes a $k$-sunflower, if $|\FF| > \lp c k^4 \rp^{m}$ for a constant $c>0$ independent of $m$ and $k$, where {\em $k$-sunflower} means a family of $k$ different sets with a common pairwise intersection. 
\end{abstract}

\maketitle

\section{Motivation, Basic Terminology and Related Facts} 
This paper verifies\footnote{ 
Extensive additional information on the proof can be found in \cite{blog}. 
} the following statement.
\begin{theorem} \label{SFC}
There exists $c \in \R_{>0}$ such that for every $k, m \in \Z_{>0}$, a family $\FF$ of sets each of cardinality at most $m$ includes a $k$-sunflower if
$
|\FF| > \lp c k^4 \rp^{m}
$. 
\qed
\end{theorem}

\noindent 
Here a $k$-sunflower means a family of $k$ different sets with a common pairwise intersection called the {\em core}.

This proves the {\em sunflower conjecture} that has been open since the sunflower lemma was shown in 1960, as particularly emphasized in \cite{E81}. 
The base $ck^4$ of the exponential lower limit is independent of the exponent $m$ as the conjecture states. Any similar lower limit must have a base at least $k-1$, since we can construct\footnote{
Such a family $\FF_{k, m}$ can be constructed as follows. The $\FF_{k, 1}$ consists of $k-1$ distinct elements of the universal set $X$. Given $\FF_{k, m-1}$ recursively, create $k-1$ new elements of $X$. For each $U \in \FF_{k, m-1}$ and new element $x$, add $U \cup \lb x \rb$ to $\FF_{k, m}$ to create the $m^{th}$ famly. It is straightforward to see it is free from a $k$-sunflower but $|\FF_{k, m}| =(k-1)^m$. 
} 
an $\FF$ without a $k$-sunflower such that $|\FF| =(k-1)^m$ for any $k \ge 2$ and $m \ge 1$.

\medskip

The rest of the section describes our basic terminology similarly to \cite{blog, KDC}. 
Let $X$ be the universal set, $n$ be its cardinality, and $\ep \in (0, 1)$ a sufficiently small positive number depending on no other variables. Denote 
\(
&& 
[b] = \lbr 1, b \rbr \cap \Z, \quad \textrm{for~} b  \in \R, 
\\ && 
{X'\choose m} = \lb U ~:~ U \subset X',~|U| =m  \rb, \quad \textrm{for~} X' \subset X 
\textrm{~and~} m \in [n], 
\\ && 
\FF^r = \underbrace{
\FF \times \FF \times \cdots \times \FF
}_r, \quad \textrm{for~} \FF \subset 2^X 
\textrm{~and~} r \in \Z_{\ge 0}, 
\\ \textrm{and} && 
\FF[S]  = \lb U ~:~ U \in \FF,~ S \subset U  \rb, \quad \textrm{for~} 
S \in 2^X. 
\)

A set means a subset of $X$, and is an {\em $m$-set} if it is in ${X \choose m}$. 
A family $\FF \subset {X \choose m}$ is said to satisfy the {\em $\Gamma(b)$-condition} if $|\FF[S]| < b^{-|S|} |\FF|$ for every nonempty set $S$. 
Denote set/family subtraction by $-$. 
For simplicity, a real interval may denote the integral interval of the same range; $e.g.$, use $\lbr 0, t \rp$ instead of $\lbr 0, t \rp \cap \Z$ if it is clear by context. 


\medskip

We have the double inequality 
\[
\lp \frac{x}{y} \rp^y \le  {x \choose y} < \lp \frac{e x}{y} \rp^y, \quad 
\textrm{for}~ x\in \Z_{>0} ~\textrm{and}~ y \in [x], 
\]
called the {\em standard estimate} of the binomial coefficient ${x \choose y}$, where $e=2.71...$ is the natural logarithm base. It is straightforward to check\footnote{
The lower bound holds since 
$
{x \choose y} = \prod_{i=0}^{y-1} \frac{x - i}{y-i}  \ge \lp \frac{x}{y} \rp^y 
$. 
The upper limit is due to $
\sqrt{2 \pi y} \lp \frac{y}{e} \rp^y
\exp \lp \frac{1}{12y+1} \rp
< y! <
\sqrt{2 \pi y} \lp \frac{y}{e} \rp^y
\exp \lp \frac{1}{12y} \rp$ for all $y \in \Z_{>0}$ shown in \cite{R55}, which implies $y! > \lp \frac{y}{e} \rp^y$ $\Rightarrow$ 
${x \choose y} \le \frac{x^y}{y!} < \lp \frac{ex}{y} \rp^y$. 
} its truth by Stirling's approximation.

\section{Proof of \refth{SFC}}  
With an independent constant $\ep \in (0, 1)$ given as above and integers $k, m \in \Z_{>2}$, denote 
\(
&& 
c = \exp \lp \ep^{-1} \rp, \quad 
h = \exp \lp c \rp, \quad 
i, j, r \in [0, m], \quad 
b_+ =(h k)^2, \quad  
\\ && 
\beta = \lc \frac{m}{c} \rc, \quad 
\gamma = (h b_+)^{-\beta-3}, \quad 
b =  k^4 \exp \lp c^2 \rp. \quad 
\) 
Assume WLOG that $\FF$ given by \refth{SFC} satisfies the $\Gamma(b)$-condition, which implies $|\FF|>b^m$.  Also $m > k^{-1} b$, otherwise the claim is clear by the sunflower lemma. 
We detect a $k$-sunflower in such an $\FF$ in the following four steps to prove the theorem: 
\be{{\hspace*{2mm}Step} 1.}
\item Show some properties of subfamilies of $\FF^2$. 
\item Define another type of objects to show a statement.
\item Select the final sunflower core $C$. 
\item Detect a $k$-sunflower in $\FF[C]$. 
\ee

{\bf Step 1.~}{\em Show some properties of subfamilies of $\FF^2$.} 
Denote 
\(
&& 
\PP_j = \lb (T, U)~:~  (T, U)\in \FF^2,~ |T \cap U| =j \rb, \quad 
\textrm{for}~j \in [0, m], 
\\ && 
\NN = \bigcup_{j=0}^\beta \PP_j, \quad 
\EE = \FF^2-  \NN, \quad 
\EE^* = \lb (T, U, C) ~:~(T, U) \in \EE,~ C \in 2^{T \cap U} \rb, 
\) 
An element of $\NN$ is called a {\em neighbor pair}.

Extend the $[\cdot]$-notation to write 
\[
\NN'[C] = \NN' \cap  \FF[C]^2, 
\qquad 
\textrm{for any $\NN' \subset \FF^2$ and $C \in 2^X$.}
\] 
Define the condition $\psi(C)$ by 
\[
\psi(C) \quad \Leftrightarrow \quad 
|\EE[C]| < b_+^{-1} \big| \NN[ C ] \big|. 
\]
See three remarks on the construction. 

\medskip 
 
\be{A)} 
\item 
$|\FF[C]| \le |\EE[C]|$ for each set $C$, 
since $\EE[C]$ includes the pair $(T, T)$ for every $T \in \FF[C]$. 
It means that if $\psi(C)$, 
\[
|\NN[C]| > 0
\quad \Rightarrow \quad 0 < |\FF[C]| \le |\EE[C]|
\quad \Rightarrow \quad |\NN[C]| > b_+. 
\]

\item $|\EE| \le |\EE^*|< b_+^{-1} \gamma^2 |\FF|^2$ meaning $|\NN| > \lp  1- b_+^{-1} \gamma^2 \rp |\FF|^2$. For, 
\[ 
|\PP_j| < {m \choose j} b^{-j} |\FF|^2 < \lp \frac{b j}{3 m} \rp^{-j} |\FF|^2 < (h^{8} k^4)^{-j} |\FF|^2, 
\] 
for all $j \in (\beta, m]$, by the $\Gamma(b)$-condition of $\FF$ and the standard estimate of ${m \choose j}$. As each $(T, U) \in \PP_j$ generates $2^j$ triples $(T, U, C) \in \EE^*$, it follows that 
\(
\frac{|\EE^*|}{|\FF|^2} &<& \sum_{j=\beta+1}^{m} \lp 2^{-1} h^{8} k^4 \rp^{-j} 
< (h^7 k^4)^{-\beta} = h^{-\beta} (h b_+)^{-2\beta}.  
\)
Further observe by $m > k^{-1} b =k^3 \exp(c^2)$ that 
\(
&& 
h^{\beta} > exp \lp \ep m \rp  > m^7 > (h b_+)^7, 
\sothat 
|\EE| \le |\EE^*|< \frac{\gamma^2}{b_+} |\FF|^2, 
\) 
justifying the claim.

\item Therefore, 
\(
\sum_{C \in 2^X,~\neg \psi(C)} ~|\NN[ C ]| 
&\le& b_+ \sum_{C \in 2^X,~\neg \psi(C)}  ~|\EE[ C ]| 
\\ &=& b_+ |\EE^*|  
< \gamma^2 |\FF|^2. 
\)
\ee

\medskip 

{\bf Step 2.~}{\em Define another type of objects to show a statement.} 
A {\em sequence of a $j$-set $S$} is a family $\lb S_0, S_1, S_2, \ldots, S_j \rb$ such that $S_0=\emptyset$ and $S_i \in {S \choose i}[S_{i-1}]$ for $i \in [j]$. 
For simplicity, use the symbol $\hat ~$ to denote a sequence of a set such as $\hat S$ of $S$. All such $\hat S$ are isomorphic to the permutations of $S$. 

Define 
\[
\sigma \lp \NN' \rp = \lb (T, U,  \hat S) ~:~ (T, U) \in \NN',~ S \in 2^{T \cap U} \rb,  
\quad \textrm{for~} \NN' \subset \NN, 
\]
including all the $(T, U, \hat S)$ incident to any $(T, U)$ and $S$, 
\[
\NN_* = \NN - \bigcup_{j \in [0, \beta],~\neg f(j)}  \PP_j, 
\eqand 
\sigma_+ (\NN') = \sigma(\NN' \cap \NN_*), 
\]
where $f(j)$ denotes the condition $|\PP_j| > \gamma |\FF|^2$.

In addition, define $\sigma_-(\NN')$ as follows. Fix a $j \in [0, \beta]$ such that $|\PP_j|$ is maximum to denote by $j_*$, and put 
\[
g(j) = 
\lb \begin{array}{cc}
j, &  \textrm{if $f(j)$ ({\em Case 1 of $j$}),} \\
j_*, &  \textrm{if $\neg f(j)$ and $j \le j_*$ ({\em Case 2}),} \\
\textrm{the largest $j' \in [0,  j-1]$ such that $f(j')$}, &  \textrm{if $\neg f(j)$ and $j > j_*$ ({\em Case 3}).} \\
\end{array} \right. 
\]
Denote by $\sigma^j_-(\NN')$ the family $\sigma\lbr \lp \NN- \NN' \rp \cap \PP_j \rbr$ if $j$ is in Case 1 or 2, and the family 
\[
\lb (T, U,  \hat S) ~:~ (T, U) \in (\NN - \NN') \cap \PP_j,~ 
S \in \bigcup_{i=0}^{g(j)} {T \cap U \choose i} \rb, 
\] 
in Case 3. 
Let 
\[
\sigma_-(\NN') = \bigcup_{j=0}^\beta \sigma^j_-(\NN'), 
\]
to complete its construction.

Set 
\[
\MM = \NN_* - \bigcup_{C \in 2^X,~\neg \psi(C)} ~\NN[ C ]. 
\] 
Then 
$
|\MM| > \lp 1 -  2 \beta \gamma \rp |\FF|^2 
$ by Remark C), and 
\beeq{eqStep2}
\MM[C] \ne \emptyset ~~\Rightarrow~~
|\NN[C]|>b_+ ~\wedge~ |\EE[C]| < b_+^{-1} |\NN[C]|, 
\quad \textrm{for any~}C \in 2^X, 
\eeq
by A).

\medskip 

We show another property of the obtained objects. 

\begin{lemma} \label{lm21}
$
|\sigma_-(\MM)| < 3^\beta \gamma  |\sigma_+(\NN)|. 
$ 
\end{lemma} 
\begin{proof}
Prove 
\beeq{eq2Step2}
\left| \sigma^j_- (\MM) \right| < 2^\beta \gamma \left| \sigma \lbr \PP_{g(j)}  \rbr \right|, 
\eeq
for each $j \in [0, \beta]$ in the three cases. 

\medskip 

{\em \noindent Case 1 of $j$:} due to $f(j)$ and Remark C), the family $\lp \NN - \MM \rp \cap \PP_j$ is smaller than 
$
\frac{\gamma^2 |\FF|^2}{\gamma |\FF|^2} =\gamma 
$ 
of $\PP_j$. Every neighbor pair $(T, U) \in \PP_j$ generates $h(j)$ triples $(T, U, \hat S) \in \sigma \lp \PP_j  \rp$, where $h(j)=\sum_{i=0}^j {j \choose i} i!$.  
As $g(j)=j$ in this case, 
\beeq{eq3Step2}
\frac{\left| \lp \NN - \MM \rp \cap \PP_j \right|}{|\PP_{g(j)}|} 
< \gamma,  
\eeq 
which implies \refeq{eq2Step2}.

\medskip 

{\em \noindent Case 2:} since $|\PP_{j_*}| > (2 \beta)^{-1} |\FF|^2$ and $j \le j_*$, 
\[
\frac{\left| \sigma^j_- (\MM) \right|}{\left| \sigma \lbr \PP_{g(j)}  \rbr \right|}=  
\frac{\left| \lp \NN - \MM \rp \cap \PP_j  \right| ~h(j)}{|\PP_{g(j)}| ~h(j_*)}  < \gamma   \quad \Rightarrow \quad \refeq{eq2Step2}.  
\]

\medskip 

{\em \noindent Case 3:} due to $|\hat S| \le g(j)$ for $(T, U, \hat S) \in \sigma_-^j(\MM)$ by construction, each $(T, U) \in (\NN - \MM) \cap \PP_j$ generates 
\[
\sum_{i=0}^{g(j)} {j \choose i} i! \le 
g(j)! \sum_{i=0}^{g(j)} {j \choose i} 
< 2^\beta g(j)! \le 2^\beta h[g(j)] 
\]
triples $(T, U, \hat S)$. Thus by \refeq{eq3Step2} holding in this case, 
\[
\frac{\left| \sigma^j_- (\MM) \right|}{\left| \sigma \lbr \PP_{g(j)}  \rbr \right|} < 
\frac{\left| \lp \NN - \MM \rp \cap \PP_j \right|~2^\beta h[g(j)]}{|\PP_{g(j)}|~h[g(j)]}< 2^\beta \gamma. 
\]

\medskip 

We have verified \refeq{eq2Step2} in every case. Noticing $f[g(j)]$ for all $j$, 
conclude that 
\(
| \sigma_-(\MM)| &=& \sum_{j=0}^\beta \left| \sigma^j_-(\MM) \right| 
< 2^\beta \gamma \sum_{j=0}^\beta \left| \sigma \lbr \PP_{g(j)} \rbr \right|  
\\ &\le&  
(\beta+1) 2^\beta \gamma \left| \sigma \lp \NN_* \rp \right|  < 3^\beta \gamma \left| \sigma_+ \lp \NN \rp \right|,  
\)
to complete the proof. 
\end{proof}

\medskip 
 
{\bf Step 3.~}{\em Select the final sunflower core $C$.} 
Extending the $[\cdot]$-notation again, write 
\[
\sigma \lp \NN' \rp [\hat D] = 
\lb (T, U,  \hat S) ~:~ (T, U, \hat S) \in \alpha(\NN'),~ \hat D \subset \hat S 
\rb,  
\]
for a sequence $\hat D$ of a set $D \in 2^X$ and $\NN' \subset \NN$. Use simlar expressions for any subfamily of $\sigma(\NN')$ such as $\sigma_+ (\NN')[\hat D]$. In addition, denote 
\[
\NN'[[D]] = \NN'[D] \cap \PP_{|D|}, \eqand 
\sigma_+ \lp \NN' \rp [[\hat D]] =  \sigma_+ \big( \NN'[[D]] \big) [\hat D]. 
\]

Find a maximal set $C$ and its sequence $\hat C$ such that 
\beeq{eqStep3}
\left| \sigma_- (\MM)[\hat C] \right| < b_+^{-\beta+|C|-1} \left| \sigma_+(\NN)[\hat C] \right|
\eeq 
There does exist such a $\hat C$, since \refeq{eqStep3} if $C=\emptyset$ by \reflm{lm21} and $|C| > \beta$ $\Rightarrow$ $\neg$\refeq{eqStep3}; the latter is true because $|C| > \beta$ $\Rightarrow$ $\NN[C]=\emptyset$ $\Rightarrow$ $\left| \sigma_+(\NN)[\hat C] \right| = \left| \sigma(\NN_*)[\hat C] \right| = 0$ meaning $\neg$\refeq{eqStep3}.

\medskip

Fix the sequence $\hat C$. The set $C$ will be the core of our $k$-sunflower detected in the last step. If $\MM[C]=\emptyset$, we would have 
\(
&& 
\sigma(\NN_*)[\hat C] = \sigma(\NN_*- \MM)[\hat C] \subset \sigma_-(\MM)[\hat C], 
\sothat 
\left| \sigma(\NN_*)[\hat C] \right| \le 
\left| \sigma_- (\MM)[\hat C] \right| < b_+^{-\beta+|C|-1} \left| \sigma_+(\NN)[\hat C] \right|
= b_+^{-\beta+|C|-1} \left| \sigma(\NN_*)[\hat C] \right|, 
\)
by \refeq{eqStep3}. Due to the contradiction and \refeq{eqStep2}, 
\beeq{eq01Step3}
|C| \le \beta, \quad 
\MM[C] \ne \emptyset, \quad 
|\NN[C]| > b_+, \eqand 
|\EE[C]| < b_+^{-1} |\NN[C]|. 
\eeq

\medskip

We verify two other properties in a lemma. 

\begin{lemma} \label{lm22} 
The set  $C$ satisfies the following two:  
\be{i)}
\item $\sum_{D \in 2^X[C] - \lb C \rb} ~|\NN_*[D]|  < 2 b_+^{-1} |\NN_*[C]|$. 
\item $|\NN[C] - \MM[C]| < 2 b_+^{-1} |\NN_*[C]|$. \qed 
\ee 
\end{lemma}

Put 
\[ 
r= |C|, \quad 
\DD_j = {X \choose j}[C],
\quad 
\hat \DD_j = 
\lb \hat D ~:~ D \in \DD_j,~\hat C \subset \hat D \rb, 
\quad \textrm{for~} j \in (r, \beta]. 
\]
to start its proof. If $r=\beta$, it means $\MM[C] = \NN_*[C]=\NN[C] \ne \emptyset$ and $\NN_*[D]=\emptyset$ holding for i) and ii).  
Thus assuming $r< \beta$, we prove the lemma by a series of three claims.

\medskip

{\noindent\bf Claim 1:} 
\[
\sum_{\hat D \in \hat \DD_j}~\left| \sigma_+ \lp \NN \rp [ \hat D ]  \right| 
< b_+^{-j+r} \left| \sigma_+ (\NN)[\hat C] \right|, \quad \textrm{for each~}  j. 
\] 
\begin{proof}
Given such a $j$, note that:  
\bdash 
\item The families $\sigma_-(\MM)[\hat D]$ are mutually disjoint over all $\hat D \in \hat \DD_j$. 

\item $\sigma_-(\MM)[\hat D] \subset \sigma_-(\MM)[\hat C]$ for each $\hat D$, since $\sigma_-(\MM)[\hat D]$ comprises the triples $(T, U, \hat S) \in \sigma_-(\MM)$ such that $\hat C \subset \hat D \subset \hat S$. 

\item 
$
\left| \sigma_+ \lp \NN \rp[\hat D] \right|  
\le b_+^{\beta - j  + 1}   \left| \sigma_-(\MM)[\hat D] \right| 
$ for each $\hat D$, 
by the maximality of $C$ such that \refeq{eqStep3}. 
\edash 

\medskip 

Thus, 
\(
\sum_{\hat D \in \hat \DD_j}~\left| \sigma_+ \lp \NN \rp[ \hat D ]  \right| 
&\le&  
\sum_{\hat D \in \hat \DD_j} b_+^{\beta - j + 1}
\left| \sigma_- (\MM) [\hat D] \right|
\\ &\le& \nonumber
b_+^{\beta  - j + 1}   \left| \sigma_- (\MM)[\hat C] \right|
\le 
b_+^{-j+r} \left| \sigma_+ (\NN)[\hat C] \right|, 
\)
proving the claim. 
\end{proof}

\medskip 

{\noindent\bf Claim 2:} 
$\left| \sigma_+ \lp \NN \rp [[\hat C]] \right| > 
\lp 1 - 3 b_+^{-1} \rp \left| \sigma_+(\NN) [\hat C] \right|$. 
\begin{proof}
Put 
\[
\sigma_A = \sigma_+ \lp \NN \rp [\hat C] - \sigma_+ \lp \NN \rp [[\hat C]], \eqand 
\sigma_B = \bigcup_{j \in (r, \beta],~\hat D \in \hat \DD_j}~\sigma_+ \lp \NN \rp [ \hat D ].  
\]
By Claim 1, 
\[
|\sigma_B| 
<  \left| \sigma_+ (\NN)[\hat C] \right| \sum_{j=r+1}^\beta b_+^{-j+r}
< \frac{1+\ep}{b_+} \left| \sigma_+ (\NN)[\hat C] \right|. 
\]

The family $\sigma_A$ comprises $(T, U, \hat S)$ for each neighbor pair $(T, U) \in \NN_*[C]$ such that $|T \cap U| > r$, set $S \in 2^{T \cap U}[C] - \lb C \rb$ and its sequence $\hat S$ such that $\hat C \subset \hat S$. If $\hat S \ne \hat C$, the triple belongs to $\sigma_B$ as well. Else, $\hat S= \hat C ~\wedge~|T \cap U| >r$, for which we can map such a $(T, U, \hat S) \in \sigma_A$ injectively to a $(T, U, \hat S') \in \sigma_B$ with $S' = T \cap U$.  By those we get $|\sigma_A| \le 2 |\sigma_B|$. 

Therefore,  
\[
|\sigma_A| \le 2 |\sigma_B|  
< \frac{2 + 2 \ep}{b_+} \left| \sigma_+ (\NN)[\hat C] \right|, 
\] 
followed by the desired inequality. 
\end{proof}

\medskip 

{\noindent\bf Claim 3:} for each $j$, 
\[
\sum_{D \in \DD_j}~|\NN_*[D]| < \frac{1+\ep}{b_+^{j-r}}  |\NN_*[C]|. 
\] 
\begin{proof} 
It follows the preceding two claims that 
\[
\sum_{\hat D \in \hat \DD_j}~\left| \sigma_+(\NN)[\hat D] \right| 
< b_+^{-j+r} \left| \sigma_+(\NN)[\hat C] \right| 
< \frac{1+\ep}{b_+^{j-r}} \left| \sigma_+(\NN)[[\hat C]] \right|. 
\]

Consider a $D \in \DD_j$ and its sequences $\hat D \in \hat \DD_j$ in its LHS. Each $(T, U) \in \NN_*[D]$ is incident to one or more triples $(T, U, \hat D) \in \bigcup_{\hat D' \in \hat \DD_j} \sigma_+ \lp \NN \rp[\hat D']$, while a $(T', U') \in \NN_*[[C]]$ to just the $(T', U', \hat C) \in \sigma_+\lp \NN  \rp[[\hat C]]$. Hence, 
\[
\frac{\sum_{D \in \DD_j}~|\NN_*[D]| }{\big| \NN_*[[C]] \big|} 
\le 
\frac{\sum_{\hat D \in \hat \DD_j}~\left| \sigma_+ \lp \NN\rp[\hat D] \right|}
{\left| \sigma_+ (\NN)[[ \hat C]] \right|} < \frac{1+\ep}{b_+^{j-r}}, 
\]
proving the claim. 
\end{proof}

\medskip 

By Claim 3, $\sum_{D \in 2^X[C] - \lb C \rb} ~|\NN_*[D]| < 2 b_+^{-1} |\NN_*[C]|$, which confirms Relation i) of the lemma.

\medskip 

We now verify ii). Due to \refeq{eqStep3}, \refeq{eq01Step3} and Claim 2, 
\beeq{eq2Step3}
\left| \sigma_-(\MM)[\hat C] \right| <  
\frac{1+\ep}{b_+} \big| \sigma_+(\NN)[[ C]] \big|, 
\eeq
meaning $\sigma_+(\NN)[[ C]] \ne \emptyset$ $\Rightarrow$ $\NN_*[[C]] \ne \emptyset$ $\Rightarrow$ $f(r)$. 

As in Claim 3, argue for \refeq{eq2Step3} that each $(T, U) \in (\NN- \MM)[C]$ is incident to one or more $(T, U, \hat S) \in \sigma_- (\MM)[\hat C]$ by construction\footnote{As defined in Step 2, 
$\sigma_- (\MM)$ includes all the triples $(T, U, \hat S) \in \sigma(\NN- \MM)$ except for those with $|\hat S| > g \lp j\rp$ for $(T, U) \in \NN- \MM$ such that $j=|T \cap U|$ is in Case 3. Since such a $j$ is greater than $g(j) \ge  r$ due to $f(r)$, the family $\sigma_- (\MM)[\hat C]$ includes the triple $(T, U, \hat C)$ for every $(T, U) \in (\NN- \MM)[C]$.  
}, while each $(T', U') \in \NN_*[[C]]$ to just the $(T', U', \hat C) \in \sigma_+(\NN)[[ \hat C ]]$. It follows that 
\(
|(\NN - \MM)[C]| &\le& \left| \sigma_-(\MM)[\hat C] \right|< \frac{1+\ep}{b_+} \big| \sigma_+(\NN)[[ C]] \big| 
\\ &=& \frac{1+\ep}{b_+} \big| \NN_*[[C]] \big| 
< \frac{2}{b_+} |\NN_*[C]|. 
\)
implying Relation ii). This completes the proof of the lemma.

\medskip 
 
{\bf Step 4.~}{\em Detect a $k$-sunflower in $\FF[C]$.} 
By \refeq{eq01Step3} and \reflm{lm22} with $\MM[C] \subset \NN_*[C] \subset \NN[C] = \FF^2[C] - \EE[C]$, we have the three properties 
\(
&& 
|\FF[C]^2 - \MM[C]| < \frac{c |\FF[C]|^2}{b_+}, \quad 
|\FF[C]|^2 > b_+, 
\\ && 
\textrm{and} \quad 
\left| \bigcup_{D \in 2^X[C]- \lb C \rb}~ \MM[D] \right| < \frac{c |\MM[C]|}{b_+}, 
\)
implying the following statement. 

\begin{proposition}
There exist  $[1- (c k)^{-1}] |\FF[C]| > 2^{-1} hk$ sets $T \in \FF[C]$ each satisfying the folloing condition.  
\bnull
\item $\phi(T)$: $T \cap U = C$ for more than $[1- (c k)^{-1}] |\FF[C]|$ sets $U \in \FF[C]$. \qed 
\enull 
\end{proposition}

\begin{figure}
\begin{small}
\be{\hspace*{5mm}}
\item Input: the set $C$ determined by Step 3. 
\item Output: a $k$-sunflower $\TT \subset \FF[C]$. 
\ee

\medskip 

\be{\qquad 1.} 
\item 
$\GG \leftarrow \FF[C]$; \hspace*{1mm}
$\TT \leftarrow \emptyset$; \hspace*{1mm}
\item \bfor $j=1$ \bto $k$ \bdo
	\be{{2}-1.}
	\item find $T_j \in \GG$ such that $\phi(T_j)$; 
	\item $\TT \leftarrow \TT \cup \lb T_j \rb$; 
	\item delete all $U \in \GG$ such that $T_j \cap U \ne C$ from $\GG$; 
	\ee 
\item \breturn$(\TT)$; 
\ee 
\end{small}
\caption{Algorithm \ltt{S}} \label{fig1}
\end{figure}

\medskip 

Noting the proposition, we construct a $k$-sunflower $\TT = \lb T_1, T_2, \ldots, T_ k \rb \subset \FF[C]$ by the algorithm \ltt{S} described in Fig.\ \ref{fig1}. Its recursive invariant is the following two conditions after Step 2-3 is executed for each $j \in [k]$: 
\be{\hspace*{1mm} a)}
\item $|\GG| \ge \lbr 1 - j (ck)^{-1} \rbr |\FF[C]|$. 
\item $T_i \cap U = C$ for all $i \in [j]$ and $U \in \GG \cup \TT - \lb T_i \rb$. 
\ee

\medskip 

a) is clearly true for every $j$, because Step 2-3 could delete $(ck)^{-1} |\FF[C]|$ or less $U$ due to $\phi(T_j)$. This also guarantees the existence of such a $T_j$ at Step 2-1 since it means $|\GG| - (ck)^{-1} |\FF[C]| > (1-\ep) |\GG|$ or more sets $T$ such that $\phi(T)$ in the current $\GG$.

As Step 2-3 eliminates all the $U \in \GG$ such $T_j \cap U \ne C$, the condition b) is true as well. 

\medskip 

These confirm that the two are true when \ltt{S} terminates with $\TT=\lb T_1, T_2, \ldots, T_k \rb$. 
As it is a $k$-sunflower in $\FF$ with the core $C$ by b), it completes our proof of \refth{SFC}.




\bibliographystyle{amsplain}

\end{document}